\documentclass[a4paper,12pt]{amsart}

\usepackage{a4wide,graphicx,amssymb}
\usepackage{tikz}
\usepackage{todonotes}

\newcommand{\PP}{\mathbb{P}}

\newif\ifdetails
\detailstrue
\newcommand{\DETAIL}[1]%
{\ifdetails\par\fbox{\begin{minipage}{0.9\linewidth}\textit{Detail:}
      #1\end{minipage}}\par\fi}
\newcommand{\TODO}[1]%
{\ifdetails\par\fbox{\begin{minipage}{0.9\linewidth}\textbf{TODO:}
      #1\end{minipage}}\par\fi}

\theoremstyle{remark}

\newtheorem{question}{Question}

\newcommand{\old}[1]{{}}

\title{Some remarks on the midrange crossing constant}

\author{\'Eva Czabarka}
\author{Inne Singgih}
\author{L\'aszl\'o Sz\'ekely}
\author{Zhiyu Wang}

\address{\'Eva Czabarka\\Department of Mathematics \\ University of South Carolina \\ Columbia SC 29212 \\ USA
\and Visiting Professor\\ Department of Pure and Applied Mathematics\\ University of Johannesburg\\
P.O. Box 524, Auckland Park, Johannesburg 2006\\South Africa}
\email{czabarka@math.sc.edu}
\address{Inne Singgih\\Department of Mathematics\\ University of South Carolina \\ Columbia SC 29212 \\ USA}
\email{isinggih@math.sc.edu}
\address{L\'aszl\'o Sz\'ekely\\Department of Mathematics \\ University of South Carolina \\ Columbia SC 29212 \\ USA
\and Visiting Professor\\ Department of Pure and Applied Mathematics\\ University of Johannesburg\\
P.O. Box 524, Auckland Park, Johannesburg 2006\\ South Africa}
\email{szekely@math.sc.edu}
\address{Zhiyu Wang\\ Department of Mathematics \\ University of South Carolina \\ Columbia SC 29212 \\ USA}
\email{zhiyuw@math.sc.edu}

\thanks{
The last three authors were  supported in part by the National Science Foundation  contract  DMS  1600811.}

\begin{document}

\begin{abstract} We verify an upper bound of Pach and T\'oth [{\it Combinatorica} \textbf{17} (1997), 427--439, {\it Discrete and Computational Geometry} \textbf{36}, (2006), 527--552] on the midrange crossing constant. Details of their $\frac{8}{9\pi^2}$ upper bound have not been
available. Our verification is different from their method and hinges on a result of Moon [{\it J. Soc.
Indust. Appl. Math.} \textbf{13}(1965), 506--510]. As Moon's result is optimal, we raise the question whether the  
midrange crossing constant is $\frac{8}{9\pi^2}$.
\end{abstract}

\maketitle

\section{Introduction}
Pach and T\'oth \cite{pachtoth} provided $n$ points in the plane and $e$ edges drawn between them under the constraints
$e/n\rightarrow \infty$ and $e=o(n^2)$,  with at most $(\frac{16}{27\pi^2}+o(1))
\frac{e^3}{n^2}$ crossings. Later they  \cite{pachtoth2} corrected the calculation for the number of crossings to  
\begin{equation} \label{truth}
\left(\frac{8}{9\pi^2}+o(1)\right)\frac{e^3}{n^2}<\left(0.0900633+o(1)\right)
\frac{e^3}{n^2}.
\end{equation} Their construction was a $\sqrt{n}\times \sqrt{n}$ grid, with the points slightly moved  into  general 
position, so that no 3 of them are collinear, 
and they joined the pairs of points with straight line segments if their distance did not exceed some number $d$. Details of neither of these calculations,
which are said to be unpleasant,
are available to the public, therefore we think that a simple alternative calculation as below is worth showing.

Let ${\rm cr}(G) $ denote the usual crossing number of the graph $G$ (for detailed definition see \cite{schaeferbook}). Let 
$\kappa(n,e)$ denote the minimum crossing number of a simple graph $G$ of order $n$ and size $e$.
According to~\cite{pachspencertoth}, there exists a positive constant $\gamma$, called the midrange crossing constant,
such that 
the limit
\begin{equation} \label{midlimit}
\lim_{n\rightarrow\infty}\kappa(n,e)\frac{n^2}{e^3}
\end{equation}
under the constraints $e/n\rightarrow \infty$ and $e=o(n^2)$, 
exists and is equal to $\gamma$.  The existence of such a constant was conjectured by Erd\H os and Guy \cite{EG}.
In fact they missed to make the second assumption \cite{pachspencertoth}. The second assumption, 
however, is needed. For completeness, we   show it next.
Note that $\kappa(n,\binom{n}{2})={\rm cr}(K_n)$.  The Harary-Hill conjecture \cite{brick,schaeferbook} implies that ${\rm cr}(K_n)=(\frac{1}{64}+o(1))n^4$.
The conjecture is supported by a corresponding construction providing the upper bound, and \cite{balogh} proves 98.5\% of the required
lower bound. Hence for $e=\binom{n}{2}$, we have 
$$\frac{1+o(1)}{8}\cdot 0.985< \kappa(n,e)\frac{n^2}{e^3},$$
contradicting (\ref{truth}) outside its range.

The first step towards proving the Erd\H os-Guy conjecture \cite{EG} was the discovery of the Crossing Lemma \cite{acns,le}. (Curiously, the papers proving the Crossing Lemma  seemed to be unaware of the Erd\H os-Guy conjecture.) The Crossing Lemma asserted that for $e\geq 4n$,
$$
\frac{1}{64}\frac{n^2}{e^3} \leq \kappa(n,e), 
$$ 
showing that $1/64\leq \gamma$. The constant $1/64$ has been improved a number of times at the cost of requiring
somewhat larger $e$. The current best constant, $1/29=0.0344...$  is due to Ackerman \cite{ackerman}.

Recently Czabarka, Reiswig, Sz\'ekely and Wang~\cite{classes} noted, that the existence  of the {midrange crossing constant} 
can be extended to the existence  of the {midrange crossing constant} $\gamma_{\mathcal C}$  for certain graph classes $\mathcal C$.
The condition is that $\mathcal C$ has to be closed for some graph operations.
Define $\kappa_\mathcal{C}(n,e)$  the minimum crossing number of a simple graph $G\in \mathcal{C}$ of order $n$ and size $e$.
The paper \cite{classes} showed that changing $\kappa(n,e)$  to $\kappa_\mathcal{C}(n,e)$  in (\ref{midlimit}), a limit under the same
condition exists,
which may or not be equal to the midrange crossing constant $\gamma$ for all graphs. For example, $\mathcal C$ can be the class of bipartite 
graphs. The existence of the midrange crossing constant for the class of bipartite graphs was needed to prove some tight
crossing number results \cite{asplund}.
Angelini, Bekos, Kaufmann, Pfister and Ueckerdt
\cite{angelini} proved a stronger version of the Crossing Lemma for bipartite graphs. Their result  implies that the midrange crossing constant for the class of bipartite graphs is at least $16/289>0.055$, making  plausible the  conjecture that the bipartite midrange crossing 
constant is bigger than the midrange crossing constant.

We utilize both the spirit and the calculations of the Moon \cite{moon} paper.
He observed that selecting $n$ points on the unit sphere independently according the uniform distribution,
and for any two points, connecting them on the shorter arc of their great circle, the expected number of crossings is $(\frac{1}{64}+o(1))n^4$,
which is asymptotically the same  as the conjectured crossing number of the complete graph in the  Harary-Hill conjecture. 
This result is truly surprising.

Our calculation uses two ideas. 
The first idea is that the construction of Pach and T\'oth is an imitation of a uniformly distributed large point set, the second is that
calculations on the sphere are simpler than calculations on the plane. We restrict the Moon construction by connecting only 
 pairs of points with distance at most $d$ for some fixed but very small $d$. This is not literally the same as the construction of Pach and T\'oth
 \cite{pachtoth}, but provides the same result. Considering that the Moon construction is optimal in expectation for $d=\pi$, one might wonder
 if it is still optimal for $d\rightarrow 0^+$.
\begin{question}
Is the midrange crossing constant $\gamma$ equal to   $\frac{8}{9\pi^2}$?
\end{question}
If the answer to this question is in the affirmative, then the rectilinear midrange crossing constant is also  $\frac{8}{9\pi^2}$. Recall that
the rectilinear crossing number is defined analogously to the crossing number, but edges have to be drawn in straight line segments
\cite{schaeferbook}. It has been known that there is a rectilinear midrange crossing constant (see the discussion in \cite{pachkplanar})
and obviously it has to be at least $\gamma$. On the other hand, as the construction of Pach and T\'oth
 \cite{pachtoth} is drawn in straight line, it forces equality if $\gamma=\frac{8}{9\pi^2}$

\section{Calculations}

Take two points $P$ and $Q$ independently from the uniform distribution on the unit sphere. The density function of the length $\alpha$
of the shorter arc connecting $P$ and $Q$ on  their great circle is $\frac{1}{2}\sin \alpha$ ($0<\alpha<\pi$). Next, select $R$ and $S$ 
as well independendently from the uniform distribution on the unit sphere. Observe that the probability of the $RS$ arc intersects the
$PQ$ arc, conditional on the length of the $PQ$ arc is $\alpha$, is
$$\frac{\alpha}{4\pi}.$$
Indeed, fixing the great circle of $P$ and $Q$, the probability that $R$ and $S$ fall into different hemispheres is $1/2$. If they fall in the
same hemisphere, then the $PQ$ and $RS$ arcs do not cross. If they fall in different hemispheres, then for any fixed $R$ and $S$,
rotating the $R$ and $S$ points around the axis connecting the poles of the great circle of $P$ and $Q$, shows that
\begin{equation} \label{alphametsz}
\PP\left[PQ \hbox{ arc crosses } RS \hbox{ arc }\big\vert\hbox{ length of }PQ=\alpha\right]=\frac{\alpha}{4\pi}.
\end{equation}
Moon \cite{moon} goes on to show from here that 
\begin{equation} \label{allmetsz}
\PP\left[PQ \hbox{ arc crosses } RS \hbox{ arc }\right]=\int_0^\pi \frac{\alpha}{4\pi}\cdot \left(\frac{1}{2}\sin \alpha\right) d\alpha =\frac{1}{8},
\end{equation}
showing that the expected number of crossings in his drawing of the complete graph is at most $\frac{1}{16}\binom{n}{2}\binom{n-2}{2}=
(\frac{1}{64}+o(1))n^4$ as he claimed.
We somewhat generalize these arguments. Consider the great circles of $P$ and $Q$, and of $R$ and $S$. With probability 1 these two great circles do not coincide, and hence have two intersection points, $T$ and $U$. Furthermore, the probability of the $PQ$ arc 
crossing the $RS$ arc does not depend on conditioning on two fixed great circles. Indeed, fixing two great circles, the length of the $PQ$ arc
as $\alpha$ and the length of the $RS$ arc as $\beta$, the probability that the $PQ$ arc crosses the $RS$ arc is
\begin{equation} \label{local}
2\cdot \frac{\alpha}{2\pi}\cdot \frac{\beta}{2\pi}.
\end{equation}

The first factor of 2 comes from deciding whether $T$ or $U$ will be the crossing point. Integrating out (\ref{local}) over arc length up to $d$, 
we obtain 
\begin{eqnarray} \label{global}
&&\PP\left[PQ \hbox{ arc crosses } RS \hbox{ arc }\hbox{ and length of }PQ\leq d \hbox{ and length of }RS\leq d\right]\\ \nonumber
& &\,\,= \int_0^d \left(\int_0^d \left(2\cdot \frac{\alpha}{2\pi}\cdot \frac{\beta}{2\pi}\right)\cdot   \frac{1}{2}\sin \alpha \,d\alpha\right)   \frac{1}{2}\sin \beta \, d\beta
= \frac{1}{8\pi^2} (\sin d- d\cos d)^2. \nonumber
\end{eqnarray}
Define now a random graph drawn on the sphere in the following way. The vertices are $n$ randomly and independently selected samples
from the uniform distribution on the unit sphere. Join vertices $P$ and $Q$ if the shorter of their great circle arc has length at most $d$,
and represent the edge between them by this arc. Based on (\ref{global}), the expected number of crossings in this drawn graphs is
\begin{eqnarray} \label{expcr}
\frac{1}{8\pi^2} \left(\sin d- d\cos d\right)^2  \cdot \frac{1}{2}\binom{n}{2}\binom{n-2}{2}    . 
\end{eqnarray}
Next we compute the expected number of edges in this graph. Recall that  the formula for the area of a cap of radius $d$ (measured on the surface) in the unit sphere is  $2\pi (1-\cos d)$. Therefore the expected number of neighbors of a vertex in our graph is
$$ 2(n-1)\pi (1-\cos d)/(4\pi),$$ and the expected number of edges in the graph  is
\begin{equation} \label{edges}
n\cdot (n-1)\cdot \frac{1-\cos d}{4}.
\end{equation}

It is not difficult to see that our random graph drawn on the sphere has size and crossing number concentrated around their
respective expected values. In fact, Moon \cite{moon} showed the concentration of the crossing number in the case $d=\pi$, 
i.e. for the complete graph. Summing up our results, (\ref{expcr}) and (\ref{edges}), for our drawing $D$ of our random graph,
we obtain
\begin{eqnarray} \label{constant}
\frac{{\rm  cr}(D)}{e(D)^3/n^2}&=&   
\frac{\frac{1}{8\pi^2} (\sin d- d\cos d)^2  \cdot \frac{1}{2}\binom{n}{2}\binom{n-2}{2}}{\Bigl(n\cdot (n-1)\cdot \frac{1-\cos d}{4} \Bigl)^3/n^2}
 \cdot (1+o(1))\\ \nonumber
&=& 
\frac{1}{\pi^2}\frac{ (\sin d- d\cos d)^2}{(1-\cos d)^3} \cdot (1+o(1)).
\end{eqnarray}
Observe that the function $\frac{ (\sin d- d\cos d)^2}{(1-\cos d)^3}$ is increasing for $0<d<\pi$. Hence, the smaller $d$ we take,
the better upper bound we have.
Taking the limit for $d\rightarrow 0^+$ of (\ref{constant}), we obtain
$$\lim_{d\rightarrow 0^+} \frac{1}{\pi^2}\frac{ (\sin d- d\cos d)^2}{(1-\cos d)^3} =\frac{8}{9\pi^2} < 0.0900633 ,$$
the correct upper bound from \cite{pachtoth2}. To formalize the graph construction, for any $\epsilon >0$, select a $d>0$ such that
$\frac{1}{\pi^2}\frac{ (\sin d- d\cos d)^2}{(1-\cos d)^3}< \frac{8}{9\pi^2}+\frac{\epsilon}{2}$. Select a large enough $n$ and $D$ 
in (\ref{constant}) such that  $\frac{{\rm  cr}(D)}{e(D)^3/n^2}<\frac{1}{\pi^2}\frac{ (\sin d- d\cos d)^2}{(1-\cos d)^3}+\frac{\epsilon}{2}$.
We are almost done---except that $D$ has a quadratic size. Take a  sufficiently bigger $N$, such that $n$ divides $N$,
and take $N/n$ copies of $D$ redrawn in the plane using stereographic projection, such that edges of different copies do not 
cross each other. Call this drawing $D'$. Clearly $\frac{{\rm  cr}(D)}{e(D)^3/n^2}=\frac{{\rm  cr}(D')}{e(D')^3/N^2}$, and the size of $D'$ satisfies the required conditions
with the appropriate choice of $N$.

\end{document}